\documentclass[12pt]{article}
\usepackage{bbm}
\usepackage{latexsym}
\usepackage{hyperref}
\usepackage{graphicx} 
\usepackage{mathrsfs}
\usepackage{amsmath}
\usepackage{graphicx}
\usepackage{booktabs}            
\usepackage{diagbox}             
\usepackage{multirow}
\usepackage{amssymb}
\usepackage{color}
\usepackage{amsthm}
\usepackage{cite}
\usepackage{indentfirst}
\usepackage{anysize}\marginsize{45mm}{45mm}{40mm}{50mm}
\usepackage{caption}
\usepackage{booktabs}
\usepackage{tabularx}
\usepackage{ragged2e} 

\usepackage{array, caption, threeparttable}

\usepackage[font=footnotesize,labelfont=bf,labelsep=none]{caption}
\newcolumntype{L}{>{\RaggedRight\hangafter=1\hangindent=0em}X}

\usepackage{enumerate}

\renewcommand\thesection{\arabic{section}.\kern -.5em}

\newtheoremstyle{lemma}{\topsep}{\topsep}%
     {}
     {}
     {\bfseries}
     {}
     {0.1em}
     {\thmname{#1}\thmnumber{ #2}\thmnote{ #3}}
\theoremstyle{lemma}  

\newtheorem{theorem}{Theorem}              
\newtheorem{lemma}[theorem]{Lemma}

\newtheorem{definition}{Definition}

\numberwithin{equation}{section}

\usepackage{latexsym, bm}
\begin{document}

\title{The vertex-pancyclicity of the simplified shuffle-cube and the vertex-bipancyclicity of the balanced shuffle-cube}

\author{Yasong Liu and Huazhong L\"{u}\thanks{Corresponding author.}\\
{\small School of Mathematical Sciences, University of Electronic Science and Technology of China,}\\
{\small Chengdu, Sichuan, 610054, PR China}\\
{\small}lvhz@uestc.edu.cn}
\date{}

\maketitle

\begin{abstract}
A graph $G$ $=$ $(V,E)$ is vertex-pancyclic if for every vertex $u$ and any integer $l$ ranging from $3$ to $|V|$, $G$ contains a cycle $C$ of length $l$ such that $u$ is on $C$. A bipartite graph $G$ $=$ $(V,E)$ is vertex-bipancyclic if for every vertex $u$ and any even integer $l$ ranging from $4$ to $|V|$, $G$ contains a cycle $C$ of length $l$ such that $u$ is on $C$. The simplified shuffle-cube and the balanced shuffle-cube, which are two variants of the shuffle-cube and are superior to shuffle-cube in terms of vertex-transitivity. In this paper, we show that the $n$-dimensional simplified shuffle-cube is vertex-pancyclic for $n\geqslant 6$, and the $n$-dimensional balanced shuffle-cube is vertex-bipancyclic for $n\geqslant 2$.

\vskip 0.1 in


\end{abstract}

\section{Introduction}
It is well known that the topology of an interconnection network, represented by a simple undirected graph, is essential in parallel systems. The hypercube $Q_{n}$ is one of the most popular and efficient interconnection networks\cite{QN}. Due to numerous appealing advantages, a lot of hypercube machines have been implemented\cite{JP}, for instance, Cosmic Cube\cite{CO}, the iPSC\cite{TH}, the Ametes S/14\cite{AM}. However, the hypercube has its own drawback, for example, its diameter is large. A number of variants of hypercube have been proposed, such as twisted cube\cite{TW}, crossed cube\cite{CR}, Lucas cube\cite{LU}, folded hypercube\cite{FO}, augmented hypercube\cite{AU}, Fibonacci cube\cite{FI}, balanced hypercube\cite{BA} and shuffle-cube\cite{ref1}.
 The shuffle-cube $SQ_{n}$ has many good properties. It is showed in\cite{ref1} the diameter of the shuffle-cube is approximately a quarter of the diameter of the hypercube of the same dimension (with the same number of vertices). Symmetry, especially vertex-transitivity, is one of the important requirements for designing high-performance interconnection networks. However, L\"{u} et al.\cite{ref8} showed that the shuffle-cube is not vertex-transitive for all $n\textgreater 2$. To address this shortcoming, in the same paper, they proposed two vertex-transitive variants of the shuffle-cube, namely, the balanced shuffle-cube and the simplified shuffle-cube. More works relating to the shuffle-cube, the balanced shuffle-cube and the simplified shuffle-cube can be found in \cite{ref8}, \cite{ref2}, \cite{ref1}, \cite{ref3}, \cite{ref4}, \cite{ref7}.
\par
A path (resp. cycle) that passes through all the vertices of a graph $G$ exactly once is called a $Hamiltonian$ $path$ (resp. $Hamiltonian$ $cycle$). For $u,v\in V(G)$, a path joining $u$ and $v$ is called a $(u,v)$-$path$ (or $(v,u)$-$path$). A graph is $Hamiltonian$ if it possesses a Hamiltonian cycle, and if there is a Hamiltonian $(u,v)$-path for any distinct vertices $u,v\in V(G)$, then it is $Hamiltonian$-$connected$. Clearly, a Hamiltonian-connected graph must be a Hamiltonian graph, but the reverse is not true, for example, Hamiltonian bipartite graphs are not Hamiltonian-connected.\par
Bondy\cite{ref9} defined that a graph $G$ is $pancyclic$ if it possesses a cycle of length $l$ (ie. an $l$-cycle) for any integer $l\in[3,|V(G)|]$. Hereafter, Randerath et al.\cite{ref10} defined a graph $G$ to be $vertex$-$pancylic$ (resp. $edge$-$pancyclic$) if for every $l\in [3,|V(G)|]$, any vertex (resp. edge) of $G$ lies on a cycle of $l$. Similarly, a bipartite graph $G$ is $bipancyclic$ if it contains  an $l$-cycle for every even integer $l\in[4,|V(G)|]$. Furthermore, a bipartite graph $G$ is $vertex$-$bipancyclic$ (resp. $edge$-$bipancyclic$) if for every even $l\in [4,|V(G)|]$, any vertex (resp. edge) of $G$ lies on a cycle of $l$. Since the embedding of cycles (and paths) is a fundamental issue in assessing the capabilities of interconnection networks \cite{AS}, it is of interest to study cycle embedding of simplified shuffle-cubes and balanced shuffle-cubes. L\"{u} et al.\cite{ref8} have showed the vertex-transitivity of both $SSQ_{n}$ and $BSQ_{n}$, and have also proved the existence of Hamiltonian cycle embeddings in both graphs. In this paper, we shall show that $SSQ_{n}$ is vertex-pancyclic for $n\geqslant 6$ and $BSQ_{n}$ is vertex-bipancyclic for $n\geqslant 2$. \par

 The rest of this paper is organized as follows. In Section 2, we give the definition and basic properties of the $n$-dimensional simplified shuffle-cube $SSQ_{n}$ and $n$-dimensional balanced shuffle-cube $BSQ_{n}$. In Section 3, we show the vertex-pancyclicity of $SSQ_{n}$. In Section 4, we prove the vertex-bipancyclicity of $BSQ_{n}$. In Section 5, we conclude the paper.\par

\vskip 0.0 in

\section{Some definitions and lemmas}
Let $G$ = $(V,E)$ be a graph, where $V$ is the vertex-set of $G$ and $E$ is the edge-set of $G$. The number of vertices of $G$ is denoted by $|G|$. A path $P$ = $\langle x_{0},x_{1},\ldots,x_{k}\rangle$ in $G$ is a sequence of distinct vertices so that there is an edge joining each pair of consecutive vertices. For convenience, we denote a path between $x_{0}$ and $x_{k}$ by $P[x_{0},x_{k}]$. If a path $C = \langle x_{0},x_{1},\ldots,x_{k}\rangle$ is such that $k \geqslant 3$, $x_{0}$ = $x_{k}$, then $C$ is said to be a $cycle$, and the length of $C$ is the number of edges contained in $C$. If there are two cycles $C_{1}$, $C_{2}$ in $G$ and $E(C_{1})\cap E(C_{2}) = e$, then $(C_{1}-e)\cup(C_{2}-e)$ is a big cycle not contain $e$.\par

Unless otherwise stated, we use $n$-bit binary sequences to label the vertices of a graph. For a vertex $u$ = $u_{n-1}u_{n-2}\ldots u_{1}u_{0}$, $u_{i} \in \{0,1\}$ for each $0 \leqslant i \leqslant n-1$, the $j$-$\mathit{prefix}$ of $u$, written by $p_{j}(u)$, is $u_{n-1}u_{n-2}\ldots u_{n-j}$, and the $k$-$\mathit{suffix}$ of $u$, written by $s_{k}(u)$, is $u_{k-1}u_{k-2}\ldots u_{1}u_{0}$. For convenience, we define the $j$-th 4-bit of a vertex $u$, denoted by $u_{4}^j$, as $u_{4}^{j} = u_{4j+1}u_{4j}u_{4j-1}u_{4j-2}$, $1\leqslant j\leqslant \frac{n-2}{4}$, and we define $u_{4}^{0} = u_{1}u_{0}$. Two different vertices $u$ and $v$, $u_{4}^{j}=v_{4}^{j}$ if they bitwise equal.\par
To recursively build simplified shuffle-cubes, we introduce a set containing 4-tuple of binary sequences, i.e. $V_{00}$ = $\{1111,0001,0010,0011\}$.\par

\vskip 0.02 in

Now we are ready to give the definition of the simplified shuffle-cube and the balanced shuffle-cube.

\begin{definition}{\bf .}\cite{ref8} The $n$-dimensional simplified shuffle-cube $SSQ_{n}$ has $2^{\frac{3n+2}{4}}$ vertices, $n$ $\equiv$ $2$(mod $4$), each of which is labeled by an $n$-bit binary string $u_{n-1}u_{n-2}\ldots u_{1}$\\
$u_{0}$. We define $SSQ_{2} \cong Q_{2}$. For $n \textgreater 2$, $SSQ_{n}$ is obtained by taking eight subcubes $SSQ_{n-4}^{i_{1}i_{2}i_{3}i_{4}}$, where $i_{1}i_{2}$ $\in$ $\{00,11\}$, $i_{3}$,$i_{4}$ $\in$ $\{0,1\}$ and all vertices of $SSQ_{n-4}^{i_{1}i_{2}i_{3}i_{4}}$ share the same $p_{4}(u)$ $=$ $i_{1}i_{2}i_{3}i_{4}$. The vertices $u$ $=$ $u_{n-1}u_{n-2}\ldots u_{1}u_{0}$ and $v$ $=$ $v_{n-1}v_{n-2}\ldots v_{1}v_{0}$ in different subcubes of dimension $n$ $-$ $4$ are adjacent in $SSQ_{n}$ iff

\begin{enumerate}
\item $s_{n-4}(u) = s_{n-4}(v)$, and\par
\item $p_{4}(u) \oplus p_{4}(v) \in V_{00}$,\\

\end{enumerate}
\end{definition}

where the notation $``\oplus"$ means bitwise addition under modulo $2$.\\

\vskip 0.0 in

\begin{figure}
\begin{minipage}[t]{0.5\linewidth}

\includegraphics[width=65mm]{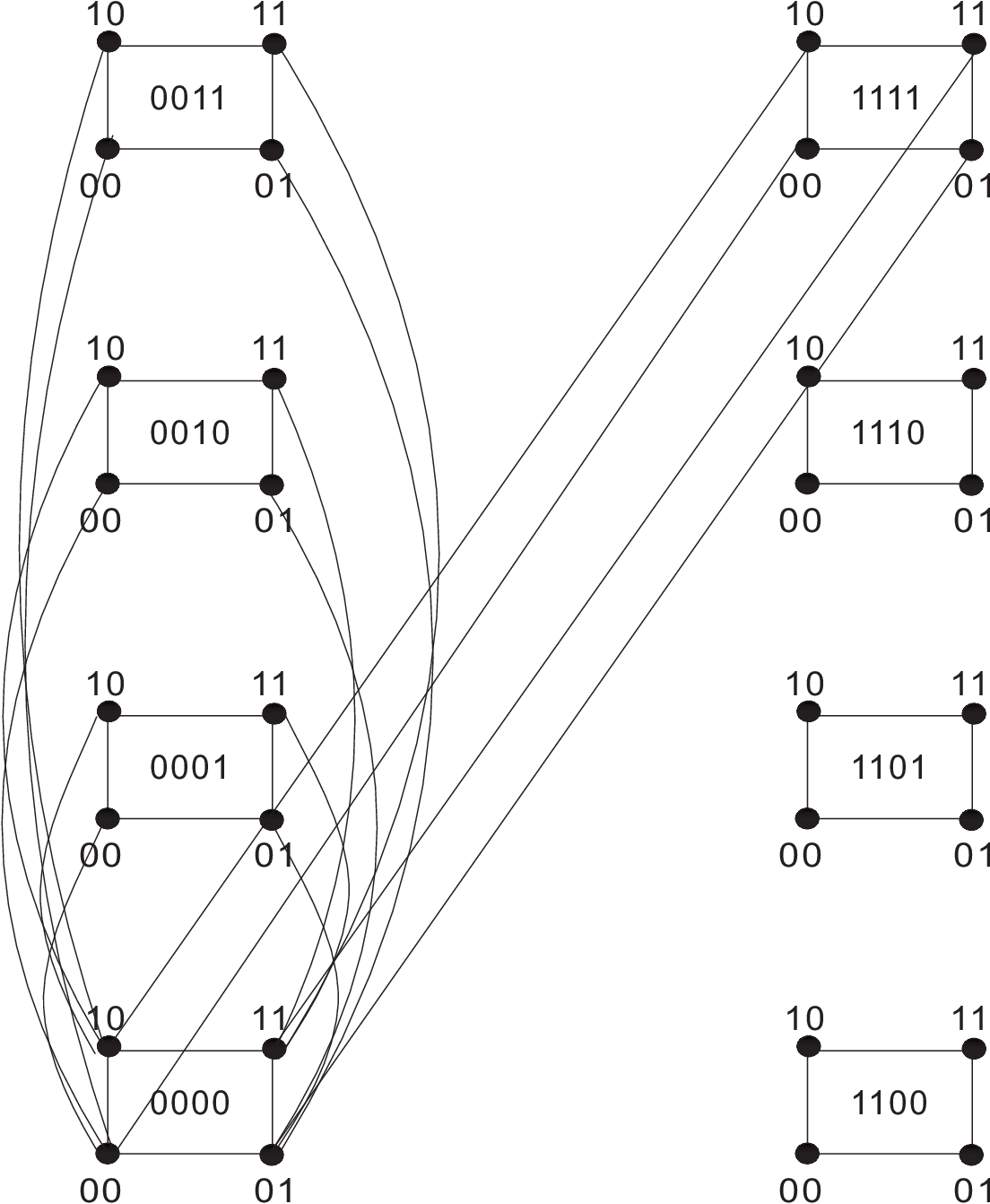}
\caption{$SSQ_{6}$.} \label{g1}
\end{minipage}
\begin{minipage}[t]{0.5\linewidth}
\centering
\includegraphics[width=70mm]{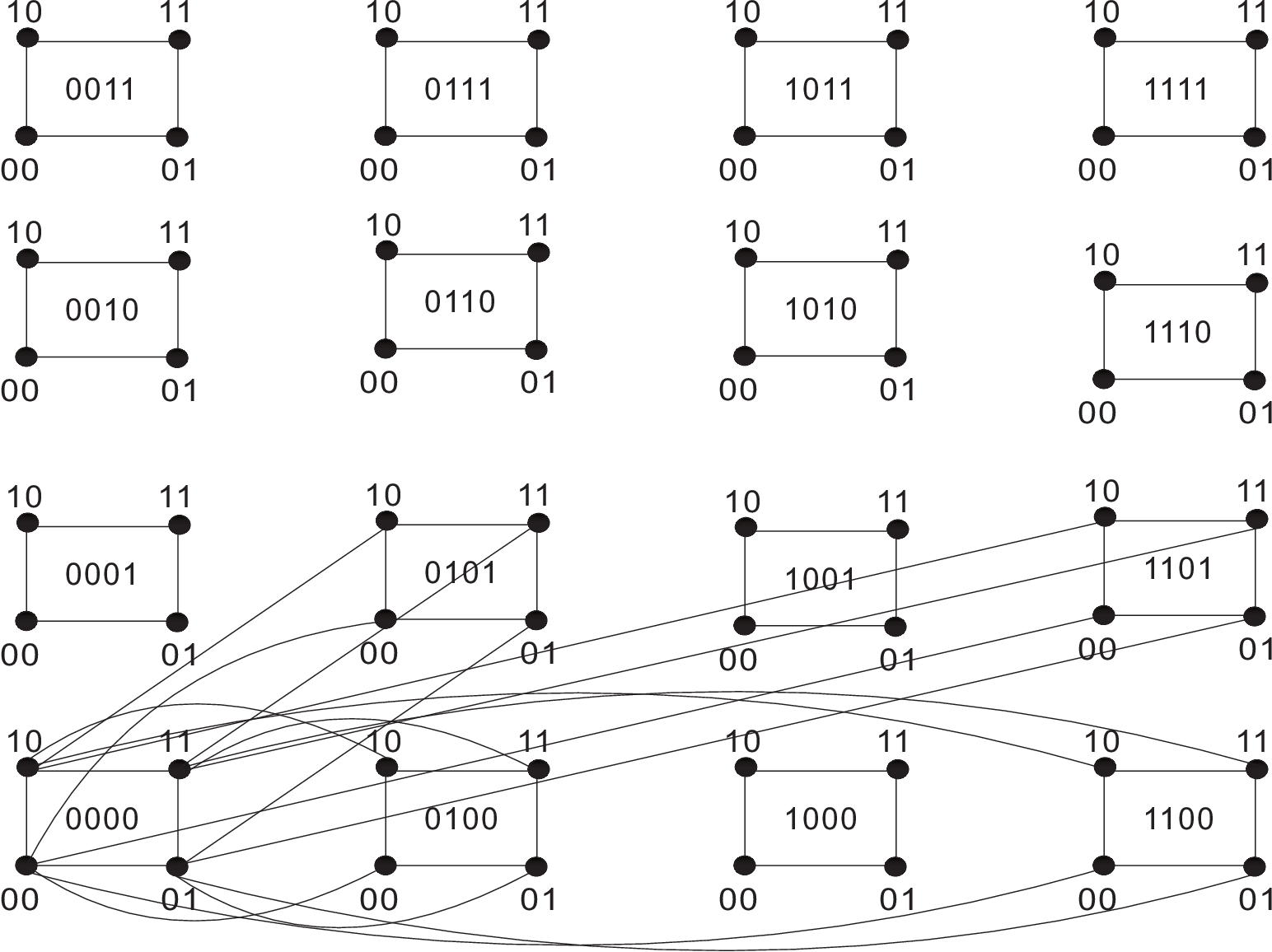}
\caption{$BSQ_{6}$.} \label{g2}
\end{minipage}
\end{figure}

By the definition above, it is clear that $SSQ_{n}$ is $n$-regular and non-bipartite. For convenience,  we let $x^{i_{1}i_{2}i_{3}i_{4}}$ be a vertex in $V(SSQ_{n-4}^{i_{1}i_{2}i_{3}i_{4}})$, and let $C^{i_{1}i_{2}i_{3}i_{4}}$ be a cycle in $SSQ_{n-4}^{i_{1}i_{2}i_{3}i_{4}}$. $SSQ_{6}$ is illustrated in Fig. \ref{g1} with only edges incident to vertices in $SSQ_{2}^{0000}$.

\begin{definition}{\bf .}\cite{ref8} The $n$-dimensional balanced shuffle-cube $BSQ_{n}$ has $2^{n}$ vertices, $n\equiv 2$(mod $4$), each of which is labeled by an $n$-bit binary string $u_{n-1}u_{n-2}\ldots u_{1}u_{0}$. We define $BSQ_{2} \cong Q_{2}$. For $n \textgreater 2$, $BSQ_{n}$ is obtained by taking sixteen subcubes $BSQ_{n-4}^{i_{1}i_{2}i_{3}i_{4}}$, where $i_{1}$,$i_{2}$,$i_{3}$,$i_{4}$ $\in$ $\{0,1\}$ and all vertices of $BSQ_{n-4}^{i_{1}i_{2}i_{3}i_{4}}$ share the same $p_{4}(u)$ $=$ $i_{1}i_{2}i_{3}i_{4}$. The vertices $u$ $=$ $u_{n-1}u_{n-2}\ldots u_{1}u_{0}$ and $v$ $=$ $v_{n-1}v_{n-2}\ldots v_{1}v_{0}$ in different subcubes of dimension $n$ $-$ $4$ are adjacent in $BSQ_{n}$ iff

\begin{enumerate}
\item $s_{n-4}(u) = s_{n-4}(v)$, and\par
\item $u_{n-2}$ and $v_{n-2}$ have different parities, and $v_{n-1}v_{n-2}=u_{n-1}u_{n-2}\pm 1$, and\par
\item $v_{n-3}v_{n-4}=u_{n-3}u_{n-4}+(-1)^{u_{n-2}}$, or $v_{n-3}v_{n-4}=u_{n-3}u_{n-4}$,\\

\end{enumerate}
\end{definition}

 where addition and substraction are under modulo 4 by regarding the two binary bits as an integer.

\vskip 0.0 in
Clearly, $BSQ_{n}$ is $n$-regular and bipartite. For convenience,  we let $x^{i_{1}i_{2}i_{3}i_{4}}$ be a vertex in $V(BSQ_{n-4}^{i_{1}i_{2}i_{3}i_{4}})$, and let $C^{i_{1}i_{2}i_{3}i_{4}}$ be a cycle in $BSQ_{n-4}^{i_{1}i_{2}i_{3}i_{4}}$. $BSQ_{6}$ is illustrated in Fig. \ref{g2} with only edges incident to vertices in $BSQ_{2}^{0000}$.
\vskip 0.0 in
The following lemmas of the $SSQ_{n}$ and the $BSQ_{n}$ will be used later.

\begin{lemma}{\bf .}\cite{ref8}
$SSQ_{n}$ is vertex-transitive whenever $n\geqslant 2$.
\end{lemma}
\begin{lemma}{\bf .}\cite{ref8}\label{SHC}
$SSQ_{n}$ is Hamiltonian-connected whenever $n\geqslant 6$.
\end{lemma}
\begin{lemma}{\bf .}\cite{ref8}
$BSQ_{n}$ is  vertex-transitive whenever $n\geqslant 2$.
\end{lemma}

\section{The vertex-pancyclicity of the simplified shuffle-cube}

We begin with the following lemma.
\begin{lemma}{\bf .}\label{cycle}
Let $uv$ be an edge of the $n$-dimensional simplified shuffle-cube $SSQ_{n}$ for $n\geqslant 6$. Then $uv$ is contained in a cycle $C$ of length 16 such that $|E(C)\cap E(SSQ_{n-4}^{i_{1}i_{2}i_{3}i_{4}})| = 1$.
\end{lemma}

\vskip 0.05 in

\noindent {\bf Proof.} {\rm Since $SSQ_{n}$ is vertex-transitive, without loss of generality, we assume that $u = 0000\ldots 0000\ldots 00$. We consider the following cases.\\
{\bf Case 1.} $v_{4}^{j} = 0001$, $0010$ or $0011$ for some integer $j$ with $1\leqslant j\leqslant \frac{n-2}{4}$, $v_{4}^{i} = 0000$ for all $i\ne j$, $v_{4}^{0} = 00$. We may assume that $v = 0000\ldots 0001\ldots 00$. Clearly, $u$ and $v$ are adjacent. It can be verified that\\
\hspace*{0.45cm}$\langle 0000\ldots 0000\ldots 00, 0000\ldots 0001\ldots 00, 1111\ldots 0001\ldots 00, 1111\ldots 0000\ldots 00,$\\
\hspace*{0.6cm}$1110\ldots 0000\ldots 00, 1110\ldots 0001\ldots 00, 1101\ldots 0001\ldots 00, 1101\ldots 0000\ldots 00,$\\
\hspace*{0.6cm}$1100\ldots 0000\ldots 00, 1100\ldots 0001\ldots 00, 0011\ldots 0001\ldots 00, 0011\ldots 0000\ldots 00,$\\
\hspace*{0.6cm}$0010\ldots 0000\ldots 00,0010\ldots 0001\ldots 00,0001\ldots 0001\ldots 00,0001\ldots 0000\ldots 00,$\\
\hspace*{0.6cm}$0000\ldots 0000\ldots 00 \rangle$\\
is the desired cycle.




\noindent {\bf Case 2.} $v_{4}^{j} = 1111$ for some integer $j$ with $1\leqslant j\leqslant \frac{n-2}{4}$, $v_{4}^{i} = 0000$ for all $i\ne j$, $v_{4}^{0} = 00$. We may assume that $v = 0000\ldots 1111\ldots 00$. It can be verified that\\
\hspace*{0.45cm}$\langle 0000\ldots 0000\ldots 00,0000\ldots 1111\ldots 00,0001\ldots 1111\ldots 00,0001\ldots 0000\ldots 00,$\\
\hspace*{0.6cm}$0010\ldots 0000\ldots 00,0010\ldots 1111\ldots 00,0011\ldots 1111\ldots 00,0011\ldots 0000\ldots 00,$\\
\hspace*{0.6cm}$1100\ldots 0000\ldots 00,1100\ldots 1111\ldots 00,1101\ldots 1111\ldots 00,1101\ldots 0000\ldots 00,$\\
\hspace*{0.6cm}$1110\ldots 0000\ldots 00,1110\ldots 1111\ldots 00,1111\ldots 1111\ldots 00,1111\ldots 0000\ldots 00,$\\
\hspace*{0.6cm}$0000\ldots 0000\ldots 00 \rangle$\\
is the desired cycle.

\noindent {\bf Case 3.} $v_{4}^{0} = 01$ or $10$, $v_{4}^{i} = 0000$ for all $i\neq 0$. We may assume that $v = 0000\ldots 0000\ldots 01$. It can be verified that\\
\hspace*{0.45cm}$\langle 0000\ldots 0000\ldots 00,0000\ldots 0000\ldots 01,0001\ldots 0000\ldots 01,0001\ldots 0000\ldots 00,$\\
\hspace*{0.6cm}$0010\ldots 0000\ldots 00,0010\ldots 0000\ldots 01,0011\ldots 0000\ldots 01,0011\ldots 0000\ldots 00,$\\
\hspace*{0.6cm}$1100\ldots 0000\ldots 00,1100\ldots 0000\ldots 01,1101\ldots 0000\ldots 01,1101\ldots 0000\ldots 00,$\\
\hspace*{0.6cm}$1110\ldots 0000\ldots 00,1110\ldots 0000\ldots 01,1111\ldots 0000\ldots 01,1111\ldots 0000\ldots 00,$\\
\hspace*{0.6cm}$0000\ldots 0000\ldots 00\rangle$\\
is the desired cycle.

The proof is complete.
 } \qed

\begin{lemma}{\bf .}\label{SSQ6}
The simplified shuffle-cube $SSQ_{6}$ is vertex-pancyclic.
\end{lemma}

\vskip 0.05 in

\noindent {\bf Proof.} {\rm Since $SSQ_{6}$ is vertex-transitive, without loss of generality, we assume that $u = 000000$. The required cycles of length from 3 to 32 passing through $u$ are listed in Table \ref{tab1}.
\begin{table*}
\setlength{\abovecaptionskip}{0cm} 
\footnotesize
\caption{\fontsize{10.5bp}{15bp}Cycles of various length passing through $u = 000000$.}
\label{tab1} 

\begin{tabularx}{1.0\linewidth}{p{1cm} L }
\toprule[1.5pt]

$C_{3}$ & $\langle 000000,001000,000100,000000 \rangle$\\
$C_{4}$ & {$\langle 000000,000001,000011,000010,000000 \rangle$} \\
$C_{5}$ & $\langle 000000,000100,000101,001001,001000,000000 \rangle$\\
$C_{6}$ & $\langle 000000,000001,000011,000111,000110,000010,000000 \rangle$\\
$C_{7}$ & $\langle 000000,000100,000101,000111,001011,001001,001000,000000 \rangle$\\
$C_{8}$ & $\langle 000000,000010,000011,000001,000101,000111,000110,000100,000000 \rangle$\\
$C_{9}$ & $\langle 000000,000100,000101,000111,000110,001010,001011,001001,001000,$$000000 \rangle$\\
$C_{10}$ & $\langle 000000,000010,000011,000001,000101,000111,001011,001010,000110,$$000100,000000\rangle$\\
$C_{11}$ & $\langle 000000,000100,000101,000111,000110,001010,001110,001111,001011,$$001001,001000,$

\ $000000 \rangle$\\
$C_{12}$ & $\langle 000000,000010,000011,000001,000101,000111,001011,001001,001000,$$001010,000110,$

$\ 000100,000000 \rangle$\\
$C_{13}$ & $\langle 000000,000100,000101,000111,000110,001010,001110,001100,001101,001111,$$001011,$

$\ 001001,001000,000000 \rangle$\\
$C_{14}$ & $\langle 000000,000100,000110,001010,001000,001100,001101,001001,001011,000111,000101,$

$\ 000001,000011,000010,000000 \rangle$\\
$C_{15}$ & $\langle 000000,000100,000101,000111,000110,001010,001110,001100,110000,110001,001101,$
$\ 001111,001011,001001,001000,000000 \rangle$\\
$C_{16}$ & $\langle 000000,000100,000110,001010,001000,001100,001110,001111,001101,$$001001,001011,$

 $\ 000111,000101,000001,000011,000010,000000 \rangle$\\
$C_{17}$ & $\langle 000000,000100,000101,000111,000110,001010,001110,001100,110000,$$110010,110011,$

$\ 110001,001101,001111,001011,001001,001000,000000 \rangle$\\
$C_{18}$ & $\langle 000000,000100,000110,001010,001000,001100,001110,110010,110011,$$001111,001101,$

$\ 001001,001011,000111,000101,000001,000011,000010,$$000000 \rangle$\\
$C_{19}$ & $\langle 000000,000100,000101,000111,000110,001010,001110,001100,110000,$$110010,110110,$

$\ 110111,110011,110001,001101,001111,001011,$$001001,$$001000,$$000000 \rangle$\\
$C_{20}$ & $\langle 000000,000100,000110,001010,001000,001100,001110,110010,110000,$$110001,110011,$

$\ 001111,001101,001001,001011,000111,000101,$$000001,$$000011,000010,000000 \rangle$\\
$C_{21}$ & $\langle 000000,000100,000101,000111,000110,001010,001110,001100,110000,110010,110110,$

$\ 110100,110101,110111,110011,110001,001101,001111,001011,001001,001000,000000 \rangle$\\
$C_{22}$ & $\langle 000000,000100,000110,001010,001000,001100,001110,110010,110000,110100,110101,$

$\ 110001,110011,001111,001101,001001,001011, 000111,000101, 000001,000011,000010,$

$\ 000000 \rangle$\\

$C_{23}$ & $\langle 000000,000100,000101,000111,000110,001010,001110,001100,110000,$$110010,110110,$

\ $110100,111000,111001,110101,110111,110011,110001,001101,001111,001011,001001,$

\ $001000,000000 \rangle$\\
$C_{24}$ & $\langle 000000,000100,000110,001010,001000,001100,001110,110010,110000,$$110100,110110,$

\ $110111,110101,110001,110011,001111,001101,001001,001011,000111,000101,000001,$

\ $000011,000010,000000 \rangle$\\
$C_{25}$ & $\langle 000000,000100,000101,000111,000110,001010,001110,001100,110000,$$110010,110110,$

\ $110100,111000,111010,111011,111001,110101,110111,110011,110001,001101,001111,$

\ $001011,001001,001000,000000 \rangle$\\
$C_{26}$ & $\langle 000000,000100,000110,001010,001000,001100,001110,110010,110000,110100,110110,$

\ $111010,111011,110111,110101,110001,110011,001111,001101,001001,001011,000111,$

\ $000101,000001,000011,$$000010,000000 \rangle$\\
$C_{27}$ & $\langle 000000,000100,000101,000111,000110,001010,001110,001100,110000,110010,110110,$

\ $110100,111000,111010,111110,111111,111011,111001,110101,110111,110011,110001,$

\ $001101,001111,001011,001001,001000,000000 \rangle$\\

        \bottomrule[1.5pt]
        \end{tabularx}
\end{table*}

\begin{table*}
\renewcommand{\thetable}{1}
\setlength{\abovecaptionskip}{0cm} 
\footnotesize
\caption{\fontsize{10.5bp}{15bp}Continued.}
\label{tab2} 

\begin{tabularx}{1.0\linewidth}{p{1cm} L }
\toprule[1.5pt]

\renewcommand{\thetable}{1}

$C_{28}$ & $\langle 000000,000100,000110,001010,001000,001100,001110,110010,110000,110100,110110,$

\ $111010,111000,111001,111011,110111,110101,110001,110011,001111,001101,001001,$

\ $001011,000111,000101,000001,000011,000010,000000 \rangle$\\
$C_{29}$ & $\langle 000000,000100,000101,000111,000110,001010,001110,001100,110000,110010,110110,$

\ $110100,111000,111010,111110,111100,111101,111111,111011,111001,110101,110111,$

\ $110011,110001,001101,001111,001011,001001,001000,000000 \rangle$\\
$C_{30}$ & $\langle 000000,000100,000110,001010,001000,001100,001110,110010,110000,110100,110110,$

\ $111010,111000,111100,111101,111001,111011,110111,110101,110001,110011,001111,$

\ $001101,001001,001011,000111,000101,000001,000011,000010,000000 \rangle$\\
$C_{31}$ & $\langle 000000,000100,000101,000111,000110,001010,001110,001100,110000,110010,110110,$

\ $110100,111000,111010,111110,111100,111101,000001,000011,111111,111011,111001,$

\ $110101,110111,110011,110001,001101,001111,001011,001001,001000,000000 \rangle$\\
$C_{32}$ & $\langle 000000,000100,000110,001010,001000,001100,001110,110010,110000,110100,110110,$

\ $111010,111000,111100,111110,111111,111101,111001,111011,110111,110101,110001,$

\ $110011,001111,001101,001001,001011,000111,000101,000001,000011,000010,000000 \rangle$\\
 \bottomrule[1.5pt]
        \end{tabularx}
\end{table*}

}  \qed

\begin{theorem}{\bf .}
The simplified shuffle-cube $SSQ_{n}$ is vertex-pancyclic whenever $n\geqslant 6$.
\end{theorem}

\vskip 0.05 in

\noindent {\bf Proof.} {\rm For convenience, let $n = 4k+2$, then $SSQ_{n}$ has $2^{3k+2}$ vertices. We prove the theorem by induction on $k$. Let $l$ be any integer with $3\leqslant l\leqslant 2^{3k+2}$, where $k\geqslant 1$. By Lemma \ref{SSQ6} we know that the theorem is true for $k = 1$. Assume that the theorem is true for $k-1$. We now consider $k$. Since $SSQ_{4k+2}$ is vertex-transitive, without loss of generality, we assume that $u = 0000\ldots0000\ldots 00$.  To complete the proof, it suffices to show that there exists a cycle of length $l$ in $SSQ_{4k+2}$ containing $u$. We consider the following cases.\\
\noindent {\bf Case 1.} $3 \leqslant l \leqslant 2^{3k-1}$. By the induction hypothesis, there exists a cycle of length $l$ in $SSQ_{n-4}^{0000}$ that contains $u$.\\
\noindent {\bf Case 2.} $2^{3k-1} + 1 \leqslant l \leqslant 2^{3k-1} + 14$. Similarly, there exists a cycle $\hat{C}$ of length $l - 14$ in $SSQ_{n-4}^{0000}$ that contains $u$. Thus, we can choose an edge $x^{0}y^{0}$ in $\hat{C}$. By Lemma \ref{cycle}, there is a cycle $C'$ of length 16 that contains $x^{0}y^{0}$ and $|E(C')\cap E(SSQ_{n-4}^{i_{1}i_{2}i_{3}i_{4}})| = 1$. Let $\hat{P} = \hat{C} - x^{0}y^{0}$ and $P' = C' -x^{0}y^{0}$. Thus $\hat{P}+P'$ is a cycle of length $l$ in $SSQ_{n}$ containing $u$.\\
\noindent {\bf Case 3.}  $2^{3k-1} + 15 \leqslant l \leqslant 2^{3k} + 12$. Clearly, there exists a cycle $\hat{C}$ of length $l - 2^{3k-1} - 12$ in $SSQ_{n-4}^{0000}$ that contains $u$. Thus, we can choose an edge $x^{0000}y^{0000}$ in $\hat{C}$. By Lemma \ref{cycle}, there is a cycle $C'=\langle x^{0000},y^{0000},x^{0001},y^{0001},x^{0010},y^{0010},x^{0011},y^{0011},\\
x^{1100},y^{1100},x^{1101},y^{1101},x^{1110},y^{1110},x^{1111},y^{1111},
x^{0000}\rangle $ of length 16 such that $x^{i_{1}i_{2}i_{3}i_{4}}\\
y^{i_{1}i_{2}i_{3}i_{4}}\in E(SSQ_{n-4}^{i_{1}i_{2}i_{3}i_{4}})$. By Theorem \ref{SHC}, there is a Hamiltonian path $P[x^{0001},y^{0001}]$ of length $2^{3k-1}-1$ in $SSQ_{n-4}^{0001}$. Therefore, we can construct a cycle of length $l$ in $SSQ_{n}$ by the following way:
$$C_{l} = (\hat{C}-x^{0000}y^{0000})\cup (C'-x^{0000}y^{0000}-x^{0001}y^{0001})\cup (P[x^{0001},y^{0001}]).$$\\
\noindent {\bf Case 4.} $2^{3k} + 13 \leqslant l \leqslant 7\times 2^{3k-1}+2$.
The proof of this case is quite analogous to that of Case 3, and the only difference is to sequentially expand the edge of $C'$ in six $SSQ_{n-4}$s (except $SSQ_{n-4}^{0000}$ and $SSQ_{n-4}^{0001}$) to a Hamiltonian path with the same ends.\\
\noindent {\bf Case 5.}  $7\times 2^{3k-1} + 3 \leqslant l \leqslant 2^{3k+2}$. By Case 1, there exists a cycle $\hat{C}$ of length $l - 7\times 2^{3k-1}$ in $SSQ_{n-4}^{0000}$ that contains $u$. Thus, we can choose an edge $x^{0000}y^{0000}$ in $\hat{C}$. Similarly, there exists a cycle $C'$ of length 16 such that $x^{i_{1}i_{2}i_{3}i_{4}}y^{i_{1}i_{2}i_{3}i_{4}}\in E(SSQ_{n-4}^{i_{1}i_{2}i_{3}i_{4}})$. By Theorem \ref{SHC}, there exists a Hamiltonian path $P[x^{i_{1}i_{2}i_{3}i_{4}},y^{i_{1}i_{2}i_{3}i_{4}}]$ of length $2^{3k-1}-1$ in $SSQ_{n-4}^{i_{1}i_{2}i_{3}i_{4}}$, where $i_{1}i_{2}i_{3}i_{4}\neq 0000$. Clearly, there exist seven Hamiltonian paths in $SSQ_{n}$ (except $SSQ_{n-4}^{0000}$). By denoting seven Hamiltonian paths as $P_{j}$, respectively, $1\leqslant j \leqslant 7$, we can construct a cycle of length $l$ in $SSQ_{n}$ in the following way:\\
$C_{l} = (\hat{C}-x^{0000}y^{0000})\cup (C'-x^{0000}y^{0000}-x^{0001}y^{0001}-x^{0010}y^{0010}-x^{0011}y^{0011}-x^{1100}y^{1100}-x^{1101}y^{1101}-x^{1110}y^{1110}-x^{1111}y^{1111})\cup (\cup_{j=1}^{7}{P_{j}})$.

Hence, $SSQ_{n}$ is vertex-pancyclic for all $n\geqslant 6$.} \qed

\section{The vertex-bipancyclicity of the balanced shuffle-cube}

\begin{lemma}{\bf .}\label{bcycle}
Let $uv$ be an edge of the $n$-dimensional balanced shuffle-cube $BSQ_{n}$ for $n\geqslant 6$. Then $uv$ is contained in a cycle $C$ of length 32 such that $|E(C)\cap E(BSQ_{n-4}^{i_{1}i_{2}i_{3}i_{4}})| = 1$, where $i_{1}$,$i_{2}$,$i_{3}$,$i_{4}$ $\in$ $\{0,1\}$.
\end{lemma}

\vskip 0.05 in

\noindent {\bf Proof.} {\rm Since $BSQ_{n}$ is vertex-transitive, without loss of generality, we assume that $u = 0000\ldots 0000\ldots 00$. We consider the following cases.\\
\noindent {\bf Case 1.} $v_{4}^{j} = 0100$ or $1100$ for some integer $j$ with $1\leqslant j\leqslant \frac{n-2}{4}$, $v_{4}^{i} = 0000$ for all $i\neq j$, $v_{4}^{0} = 00$. We may assume that $v = 0000\ldots 0100\ldots 00$. It can be verified that\\
\noindent\hspace*{0.45cm}$\langle 0000\ldots 0000\ldots 00,0000\ldots 0100\ldots 00,1101\ldots 0100\ldots 00,1101\ldots 0000\ldots 00,$\\
\hspace*{0.6cm}$1001\ldots 0000\ldots 00,1001\ldots 0100\ldots 00,0101\ldots 0100\ldots 00,0101\ldots 0000\ldots 00,$\\
\hspace*{0.6cm}$0001\ldots 0000\ldots 00,0001\ldots 0100\ldots 00,1110\ldots 0100\ldots 00,1110\ldots 0000\ldots 00,$\\
\hspace*{0.6cm}$1010\ldots 0000\ldots 00,1010\ldots 0100\ldots 00,0110\ldots 0100\ldots 00,0110\ldots 0000\ldots 00,$\\
\hspace*{0.6cm}$0010\ldots 0000\ldots 00,0010\ldots 0100\ldots 00,1111\ldots 0100\ldots 00,1111\ldots 0000\ldots 00,$\\
\hspace*{0.6cm}$1011\ldots 0000\ldots 00,1011\ldots 0100\ldots 00,0111\ldots 0100\ldots 00,0111\ldots 0000\ldots 00,$\\
\hspace*{0.6cm}$0011\ldots 0000\ldots 00,0011\ldots 0100\ldots 00,1100\ldots 0100\ldots 00,1100\ldots 0000\ldots00,$\\
\hspace*{0.6cm}$1000\ldots 0000\ldots 00,1000\ldots 0100\ldots 00,0100\ldots 0100\ldots00,0100\ldots 0000\ldots 00,$\\
\hspace*{0.6cm}$0000\ldots 0000\ldots 00\rangle$\\
\noindent is the desired cycle.












\noindent {\bf Case 2.} $v_{4}^{j} = 0101$ or $1101$ for some integer $j$ with $1\leqslant j\leqslant \frac{n-2}{4}$, $v_{4}^{i} = 0000$ for all $i\neq j$, $v_{4}^{0} = 00$. We may assume that $v = 0000\ldots 0101\ldots 00$. It can be verified that\\
\noindent\hspace*{0.45cm}$\langle 0000\ldots 0000\ldots 00,0000\ldots 0101\ldots 00,0100\ldots 0101\ldots 00,0100\ldots 0000\ldots 00,$\\
\hspace*{0.6cm}$1000\ldots 0000\ldots 00,1000\ldots 0101\ldots 00,1100\ldots 0101\ldots 00,1100\ldots 0000\ldots 00,$\\
\hspace*{0.6cm}$1011\ldots 0000\ldots 00,1011\ldots 0101\ldots 00,1111\ldots 0101\ldots 00,1111\ldots 0000\ldots 00,$\\
\hspace*{0.6cm}$0011\ldots 0000\ldots 00,0011\ldots 0101\ldots 00,0111\ldots 0101\ldots 00,0111\ldots 0000\ldots 00,$\\
\hspace*{0.6cm}$0010\ldots 0000\ldots 00,0010\ldots 0101\ldots 00,1110\ldots 0101\ldots 00,1110\ldots 0000\ldots 00,$\\
\hspace*{0.6cm}$1010\ldots 0000\ldots 00,1010\ldots 0101\ldots 00,0110\ldots 0101\ldots 00,0110\ldots 0000\ldots00,$\\
\hspace*{0.6cm}$0001\ldots 0000\ldots 00,0001\ldots 0101\ldots 00,1101\ldots 0101\ldots 00,1101\ldots 0000\ldots 00,$\\
\hspace*{0.6cm}$1001\ldots 0000\ldots 00,1001\ldots 0101\ldots 00,0101\ldots 0101\ldots 00,0101\ldots 0000\ldots 00,$\\
\hspace*{0.6cm}$0000\ldots 0000\ldots 00 \rangle$\\
\noindent is the desired cycle.











\noindent {\bf Case 3.} $v_{4}^{0} = 01$ or $10$, $v_{4}^{i} = 0000$ for all $i\neq 0$. We may assume that $v=0000\ldots 0000\ldots 01$. It can be verified that\\
\noindent$\hspace*{0.45cm}\langle 0000\ldots 0000\ldots 00,0000\ldots 0000\ldots 01,0100\ldots 0000\ldots 01,0100\ldots 0000\ldots 00,$\\
\hspace*{0.6cm}$1000\ldots 0000\ldots 00,1000\ldots 0000\ldots 01,1100\ldots 0000\ldots 01,1100\ldots 0000\ldots 00,$\\
\hspace*{0.6cm}$0011\ldots 0000\ldots 00,0011\ldots 0000\ldots 01,0111\ldots 0000\ldots 01,0111\ldots 0000\ldots 00,$\\
\hspace*{0.6cm}$1011\ldots 0000\ldots00,1011\ldots 0000\ldots 01,1111\ldots 0000\ldots 01,1111\ldots 0000\ldots 00,$\\
\hspace*{0.6cm}$0010\ldots 0000\ldots 00,0010\ldots 0000\ldots 01,0110\ldots 0000\ldots 01,0110\ldots 0000\ldots 00,$\\
\hspace*{0.6cm}$1010\ldots 0000\ldots 00,1010\ldots 0000\ldots 01,1110\ldots 0000\ldots 01,1110\ldots 0000\ldots 00,$\\
\hspace*{0.6cm}$0001\ldots 0000\ldots 00,0001\ldots 0000\ldots 01,0101\ldots 0000\ldots 01,0101\ldots 0000\ldots 00,$\\
\hspace*{0.6cm}$1001\ldots 0000\ldots 00,1001\ldots 0000\ldots 01,1101\ldots 0000\ldots 01,1101\ldots 0000\ldots 00,$\\
\hspace*{0.6cm}$0000\ldots 0000\ldots 00\rangle $\\
\noindent is the desired cycle.

The proof is complete.
}  \qed

\begin{lemma}{\bf .}\label{BHC}
For any edge $e \in E(BSQ_{n})$, there exists a Hamiltonian cycle $C$ in $BSQ_{n}$ containing $e$, $n\geqslant2$.
\end{lemma}

\vskip 0.05 in

\noindent {\bf Proof.} {\rm For convenience, let $n = 4k+2$. We prove the lemma by induction on $k$. Clearly, the lemma is true for $BSQ_{2}$. Assume that the lemma is true for $k - 1$. We now consider $k$. Let $e$ be an edge in $BSQ_{4k+2}$. It suffices to show that there exists a Hamiltonian cycle in $BSQ_{4k+2}$ containing $e$. We consider the following two cases.\\
\noindent {\bf Case 1.} $e\in E(BSQ_{4k-2}^{i_{1}i_{2}i_{3}i_{4}})$ for any $i_{1}$,$i_{2}$,$i_{3}$,$i_{4}$ $\in$ $\{0,1\}$. Without loss of generality, we assume that $e\in E(BSQ_{4k-2}^{0000})$. By the induction hypothesis, there exists a Hamiltonian cycle $C_{0}$ in $BSQ_{4k-2}^{0000}$ containing $e$. We choose another edge $e_{0}\in C_{0}$, by Lemma \ref{bcycle}, $e_{0}$ is contained in a cycle $C$ of length 32 such that $|E(C)\cap E(BSQ_{4k-2}^{i_{1}i_{2}i_{3}i_{4}})| = 1$, where $i_{1}$,$i_{2}$,$i_{3}$,$i_{4}$ $\in$ $\{0,1\}$. Similarly, we can find 31 Hamiltonian cycles $C_{i}$ in $BSQ_{4k+2}$, where $1\leqslant i\leqslant 31$ (except $BSQ_{4k-2}^{0000}$). For convenience, we may assume that $e_{i}\in E(C_{i})$ and $e_{i}\in E(C)$ when $0\leqslant i \leqslant 31$. Then the Hamiltonian cycle $\hat{C}$ in $BSQ_{4k+2}$ can be constructed as follows:
$$\hat{C} = (\cup_{i=0}^{31}{(C_{i}-e_{i})})\cup (C-\cup_{i=0}^{31}{e_{i}}).$$
\noindent {\bf Case 2.} $e\notin E(BSQ_{4k-2}^{i_{1}i_{2}i_{3}i_{4}})$. By Lemma \ref{bcycle}, $e$ is contained in a cycle $C$ of length 32 such that $|E(C)\cap E(BSQ_{4k-2}^{i_{1}i_{2}i_{3}i_{4}})| = 1$, where $i_{1}$,$i_{2}$,$i_{3}$,$i_{4}$ $\in$ $\{0,1\}$. According to the induction hypothesis, we can find 32 Hamiltonian cycles $C_{i}$ in  $BSQ_{4k+2}$, where $0\leqslant i\leqslant 31$. For convenience, we may assume that $e_{i}\in E(C_{i})$ and $e_{i}\in E(C)$. Then the Hamiltonian cycle $\hat{C}$ in $BSQ_{4k+2}$ can be constructed as follows:
$$\hat{C} = (\cup_{i=0}^{31}{(C_{i}-e_{i})})\cup (C-\cup_{i=0}^{31}{e_{i}}).$$

The proof is complete.
}  \qed

\begin{lemma}{\bf .}\label{BSQ6}
The  balanced shuffle-cube $BSQ_{6}$ is vertex-bipancyclic.
\end{lemma}

\vskip 0.05 in

\noindent {\bf Proof.} {\rm Since $BSQ_{6}$ is vertex-transitive, without loss of generality, we assume that $u = 000000$. The required cycles of length from 4 to 64 passing through $u$ are listed in Table \ref{tab3}.  \qed

\begin{table*}
\renewcommand{\thetable}{2}
\setlength{\abovecaptionskip}{0cm} 
\footnotesize
\caption{\fontsize{10.5bp}{15bp}Cycles of various length passing through $u = 000000$.}
\label{tab3}
\begin{tabularx}{1.0\linewidth}{p{1cm} L }
\toprule[1.5pt]
$C_{4}$ & $\langle 000000,000001,000011,000010,000000 \rangle$\\
$C_{6}$ & $\langle 000000,000001,000011,010111,010110,000010,000000 \rangle$\\
$C_{8}$ & $\langle 000000,000001,000011,010111,010101,010100,010110,000010,000000 \rangle$\\
$C_{10}$ & $\langle 000000,000001,000011,010111,010101,100101,100100,010100,010110,000010,000000 \rangle$\\
$C_{12}$ & $\langle 000000,000001,000011,010111,010101,100101,100111,100110,100100,010100,010110,$

\ $000010,000000 \rangle$\\
$C_{14}$ & $\langle 000000,000001,000011,010111,010101,100101,100111,110111,110110,100110,100100,$

\ $010100,010110,000010,000000 \rangle$\\
$C_{16}$ & $\langle 000000,000001,000011,010111,010101,100101,100111,110111,110101,110100,110110,$

\ $100110,100100,010100,010110,000010,000000 \rangle$\\
$C_{18}$ & $\langle 000000,000001,000011,010011,010001,100001,100011,110011,110001,001101,001100,$

\ $110000,110010,100010,100000,010000,010010,000010,000000 \rangle$\\
$C_{20}$ & $\langle 000000,000001,000011,010011,010001,100001,100011,110011,110001,001101,001111,$

\ $001110,001100,110000,110010,100010,100000,010000,010010,000010,000000 \rangle$\\
$C_{22}$ & $\langle 000000,000001,000011,010011,010001,100001,100011,110011,110001,001101,001111,$

\ $011111,011110,001110,001100,110000,110010,100010,100000,010000,010010,000010,$

\ $000000 \rangle$\\
$C_{24}$ & $\langle 000000,000001,000011,010011,010001,100001,100011,110011,110001,001101,001111,$

\ $011111,011101,011100,011110,001110,001100,110000,110010,100010,100000,010000,$

\ $010010,000010,000000\rangle$\\
$C_{26}$ & $\langle 000000,000001,000011,010011,010001,100001,100011,110011,110001,001101,001111,$

\ $011111,011101,101101,101100,011100,011110,001110,001100,110000,110010,100010,$

\ $100000,010000,010010,000010,000000\rangle$\\
$C_{28}$ & $\langle 000000,000001,000011,010011,010001,100001,100011,110011,110001,001101,001111,$

\ $011111,011101,101101,101111,101110,101100,011100,011110,001110,001100,110000,$

\ $110010,100010,100000,010000,010010,000010,000000 \rangle$\\
$C_{30}$ & $\langle 000000,000001,000011,010011,010001,100001,100011,110011,110001,001101,001111,$

\ $011111,011101,101101,101111,111111,111110,101110,101100,011100,011110,001110,$

\ $001100,110000,110010,100010,100000,010000,010010,000010,000000 \rangle$\\
$C_{32}$ & $\langle 000000,000001,000011,010011,010001,100001,100011,110011,110001,001101,001111,$

\ $011111,011101,101101,101111,111111,111101,111100,111110,101110,101100,011100,$

\ $011110,001110,001100,110000,110010,100010,100000,010000,010010,000010,000000 \rangle$\\
$C_{34}$ & $\langle 000000,000001,000011,010011,010001,100001,100011,110011,110001,001101,001111,$

\ $011111,011101,101101,101111,111111,111101,001001,001000,111100,111110,101110,$

\ $101100,011100,011110,001110,001100,110000,110010,100010,100000,010000,010010,$

\ $000010,000000 \rangle$\\
$C_{36}$ & $\langle 000000,000001,000011,010011,010001,100001,100011,110011,110001,001101,001111,$

\ $011111,011101,101101,101111,111111,111101,001001,001011,001010,001000,111100,$

\ $111110,101110,101100,011100,011110,001110,001100,110000,110010,100010,100000,$

\ $010000,010010,000010,000000 \rangle$\\
$C_{38}$ & $\langle 000000,000001,000011,010011,010001,100001,100011,110011,110001,001101,001111,$

\ $011111,011101,101101,101111,111111,111101,001001,001011,011011,011010,001010,$

\ $001000,111100,111110,101110,101100,011100,011110,001110,001100,110000,110010,$

\ $100010,100000,010000,010010,000010,000000 \rangle$\\

        \bottomrule[1.5pt]
        \end{tabularx}
\end{table*}

\begin{table*}
\renewcommand{\thetable}{2}
\setlength{\abovecaptionskip}{0cm} 
\footnotesize
\caption{\fontsize{10.5bp}{15bp}Continued.}
\label{tab4}
\begin{tabularx}{1.0\linewidth}{p{1cm} L }
\toprule[1.5pt]
$C_{40}$ & $\langle 000000,000001,000011,010011,010001,100001,100011,110011,110001,001101,001111,$

\ $011111,011101,101101,101111,111111,111101,001001,001011,011011,011001,011000,$

\ $011010,001010,001000,111100,111110,101110,101100,011100,011110,001110,001100,$

\ $110000,110010,100010,100000,010000,010010,000010,000000 \rangle$\\
$C_{42}$ & $\langle 000000,000001,000011,010011,010001,100001,100011,110011,110001,001101,001111,$

\ $011111,011101,101101,101111,111111,111101,001001,001011,011011,011001,101001,$

\ $101000,011000,011010,001010,001000,111100,111110,101110,101100,011100,011110,$

\ $001110,001100,110000,110010,100010,100000,010000,010010,000010,000000 \rangle$\\
$C_{44}$ & $\langle 000000,000001,000011,010011,010001,100001,100011,110011,110001,001101,001111,$

\ $011111,011101,101101,101111,111111,111101,001001,001011,011011,011001,101001,$

\ $101011,101010,101000,011000,011010,001010,001000,111100,111110,101110,101100,$

\ $011100,011110,001110,001100,110000,110010,100010,100000,010000,010010,000010,$

\ $000000 \rangle$\\
$C_{46}$ & $\langle 000000,000001,000011,010011,010001,100001,100011,110011,110001,001101,001111,$

\ $011111,011101,101101,101111,111111,111101,001001,001011,011011,011001,101001,$

\ $101011,111011,111010,101010,101000,011000,011010,001010,001000,111100,111110,$

\ $101110,101100,011100,011110,001110,001100,110000,110010,100010,100000,010000,$

\ $010010,000010,000000 \rangle$\\
$C_{48}$ & $\langle 000000,000001,000011,010011,010011,100001,100011,110011,110001,001101,001111,$

\ $011111,011101,101101,101111,111111,111101,001001,001011,011011,011001,101001,$

\ $101011,111011,111001,111000,111010,101010,101000,011000,011010,001010,001000,$

\ $111100,111110,101110,101100,011100,011110,001110,001100,110000,110010,100010,$

\ $100000,010000,010010,000010,000000 \rangle$\\
$C_{50}$ & $\langle 000000,000001,000011,010011,010001,100001,100011,110011,110001,001101,001111,$

\ $011111,011101,101101,101111,111111,111101,001001,001011,011011,011001,101001,$

\ $101011,111011,111001,000101,000100,111000,111010,101010,101000,011000,011010,$

\ $001010,001000,111100,111110,101110,101100,011100,011110,001110,001100,110000,$

\ $110010,100010,100000,010000,010010,000010,000000 \rangle$\\
$C_{52}$ & $\langle 000000,000001,000011,010011,010001,100001,100011,110011,110001,001101,001111,$

\ $011111,011101,101101,101111,111111,111101,001001,001011,011011,011001,101001,$

\ $101011,111011,111001,000101,000111,000110,000100,111000,111010,101010,101000,$

\ $011000,011010,001010,001000,111100,111110,101110,101100,011100,011110,001110,$

\ $001100,100000,100010,100010,100000,010000,010010,000010,000000 \rangle$\\
$C_{54}$ & $\langle 000000,000001,000011,010011,011101,100001,100011,110011,110001,001101,001111,$

\ $011111,011101,101101,101111,111111,111101,001001,001011,011011,011001,101001,$

\ $101011,111011,111001,000101,000111,010111,010110,000110,000100,111000,111010,$

\ $101010,101000,011000,011010,001010,001000,111100,111110,101110,101100,011100,$

\ $011110,001110,001100,110000,110010,100010,100000,010000,010010,000010,000000 \rangle$\\
        \bottomrule[1.5pt]
        \end{tabularx}
\end{table*}

\vskip 0.05 in

\begin{theorem}{\bf .}
The  balanced shuffle-cube $BSQ_{n}$ is vertex-bipancyclic where $n\geqslant 2$.
\end{theorem}

\begin{table*}
\renewcommand{\thetable}{2}
\setlength{\abovecaptionskip}{0cm} 
\footnotesize
\caption{\fontsize{10.5bp}{15bp}Continued.}
\label{tab5}
\begin{tabularx}{1.0\linewidth}{p{1cm} L }
\toprule[1.5pt]
$C_{56}$ & $\langle 000000,000001,000011,010011,010001,100001,100011,110011,110001,001101,001111,$

\ $011111,011101,101101,101111,111111,111101,001001,001011,011011,011001,101001,$

\ $101011,111011,111001,000101,000111,010111,010101,010100,010110,000110,000100,$

\ $111000,111010,101010,101000,011000,011010,001010,001000,111100,111110,101110,$

\ $101100,011100,011110,001110,001100,110000,110010,100010,100000,010000,010010,$

\ $000010,000000 \rangle$\\
$C_{58}$ & $\langle 000000,000001,000011,010011,010001,100001,100011,110011,110001,001101,001111,$

\ $011111,011101,101101,101111,111111,111101,001001,001011,011011,011001,101001,$

\ $101011,111011,111001,000101,000111,010111,010101,100101,100100,010100,010110,$

\ $000110,000100,111000,111010,101010,101000,011000,011010,001010,001000,111100,$

\ $111110,101110,101100,011100,011110,001110,001100,110000,110010,100010,100000,$

\ $010000,010010,000010,000000 \rangle$\\
$C_{60}$ & $\langle 000000,000001,000011,010011,010001,100001,100011,110011,110001,001101,001111,$

\ $011111,011101,101101,101111,111111,111101,001001,001011,011011,011001,101001,$

\ $101011,111011,111001,000101,000111,010111,010101,100101,100111,100110,100100,$

\ $010100,010110,000110,000100,111000,111010,101010,101000,011000,011010,001010,$

\ $001000,111100,111110,101110,101100,011100,011110,001110,001100,110000,110010,$

\ $100010,100000,010000,010010,000010,000000 \rangle$\\
$C_{62}$ & $\langle 000000,000001,000011,010011,010001,100001,100011,110011,110001,001101,001111,$

\ $011111,011101,101101,101111,111111,111101,001001,001011,011011,011001,101001,$

\ $101011,111011,111001,000101,000111,010111,010101,100101,100111,110111,110110,$

\ $100110,100100,010100,010110,000110,000100,111000,111010,101010,101000,011000,$

\ $011010,001010,001000,111100,111110,101110,101100,011100,011110,001110,001100,$

\ $110000,110010,100010,100000,010000,010010,000010,000000 \rangle$\\
$C_{64}$ & $\langle 000000,000001,000011,010011,010001,100001,100011,110011,110001,001101,001111,$

\ $011111,011101,101101,101111,111111,111101,001001,001011,011011,011001,101001,$

\ $101011,111011,111001,000101,000111,010111,010101,100101,100111,110111,110101,$

\ $110100,110110,100110,100100,010100,010110,000110,000100,111000,111010,101010,$

\ $101000,011000,011010,001010,001000,111100,111110,101110,101100,011100,011110,$

\ $001110,001100,110000,110010,100010,100000,010000,010010,000010,000000 \rangle$\\
 \bottomrule[1.5pt]
        \end{tabularx}
\end{table*}
}

\vskip 0.05 in

\noindent {\bf Proof.} {\rm For convenience, let $n = 4k+2$. We prove the theorem by induction on $k$. Let $l$ be any even integer with $4\leqslant l\leqslant 2^{4k+2}$, where $k\geqslant 0$. Clearly, the theorem is true for $k=0$. By Lemma \ref{BSQ6}, we know that the theorem is true for $k = 1$. Assume that the lemma is true for $k-1$. We now consider $k$. Since $BSQ_{4k+2}$ is vertex-transitive. Without loss of generality, we assume that $u = 0000\ldots 0000\ldots00$. It suffices to show that there exists a cycle of length $l$ in $BSQ_{n}$ containing $u$. We consider the following cases.\\
{\bf Case 1.}  $4 \leqslant l \leqslant 2^{4k-2}$. By the induction hypothesis, there exists a cycle of length $l$ in $BSQ_{n-4}^{0000}$ that contains $u$.\\
{\bf Case 2.} $2^{4k-2} + 2 \leqslant l \leqslant 2^{4k-2} + 30$. Similarly, there exists a cycle $\hat{C}$ of length $l - 30$ in $BSQ_{n-4}^{0000}$ that contains $u$. Thus, we can choose an edge $x^0y^0$ in $\hat{C}$. By Lemma \ref{bcycle}, there exists a cycle $C'$ of length 32 such that $|E(C')\cap E(BSQ_{n-4}^{i_{1}i_{2}i_{3}i_{4}})| = 1$. Let $\hat{P} = \hat{C} - x^{0}y^{0}$ and $P' = C' -x^{0}y^{0}$. Thus $\hat{P}+P'$ is a cycle of length $l$ in $BSQ_{n}$ containing $u$.\\
{\bf Case 3.}  $2^{4k-2} + 32 \leqslant l \leqslant 2^{4k-1} + 28$. Again, there exists a cycle $\hat{C}$ of length $l - 2^{4k-2} - 28$ in $BSQ_{n-4}^{0000}$ that contains $u$. Thus, we can choose an edge $x^{0000}y^{0000}$ in $\hat{C}$. By Lemma \ref{bcycle}, there exists a cycle $C'=\langle x^{0000},y^{0000},x^{0100}, y^{0100},x^{1000},y^{1000},x^{1100},\\
y^{1100},x^{0011},y^{0011},x^{0111},y^{0111},x^{1011},y^{1011},x^{1111},y^{1111},x^{0010},y^{0010},x^{0110},y^{0110},x^{1010},\\
y^{1010},x^{1110},y^{1110},x^{0001},y^{0001},x^{0101},y^{0101},x^{1001},y^{1001},x^{1101},y^{1101},x^{0000}\rangle $ of length 32 \\
such that $x^{i_{1}i_{2}i_{3}i_{4}}y^{i_{1}i_{2}i_{3}i_{4}}\in E(BSQ_{n-4}^{i_{1}i_{2}i_{3}i_{4}})$. By Theorem \ref{BHC}, there exists a Hamiltonian cycle $C^{0100}$ of length $2^{4k-2}$ in $BSQ_{n-4}^{0100}$ contains $x^{0100}y^{0100}$. Therefore, we can construct a cycle of length $l$ in $BSQ_{n}$ as follows:\\
$$C_{l} = (\hat{C}-x^{0000}y^{0000})\cup (C'-x^{0000}y^{0000}-x^{0100}y^{0100})\cup (C^{0100}-x^{0100}y^{0100}).$$\\
{\bf Case 4.} $2^{4k-1} + 30 \leqslant l \leqslant 15\times 2^{4k-2}+2$.
The proof of this case is quite analogous to that of Case 3, the only difference is to sequentially expand the edge of $C'$ in thirty $BSQ_{n-4}$s (except $BSQ_{n-4}^{0000}$ and $BSQ_{n-4}^{0100}$) to a Hamiltonian path with the same ends.\\
{\bf Case 5.} $15\times 2^{4k-2} + 4 \leqslant l \leqslant 2^{4k+2}$. There exists a cycle $\hat{C}$ of length $l - 15\times 2^{4k-2}$ in $BSQ_{n-4}^{0000}$ that contains $u$. Thus, we can choose an edge $x^{0000}y^{0000}$ in $\hat{C}$. Similarly, there exists a cycle $C'$ of length 32 such that $x^{i_{1}i_{2}i_{3}i_{4}}y^{i_{1}i_{2}i_{3}i_{4}}\in E(BSQ_{n-4}^{i_{1}i_{2}i_{3}i_{4}})$. By Theorem \ref{BHC}, there exist 31 Hamiltonian cycles $C^{i_{1}i_{2}i_{3}i_{4}}$ of length $2^{4k-2}$ in $BSQ_{n}$ (except $BSQ_{n-4}^{0000}$). By denoting 31 Hamiltonian cycles as $C_{j}$ and $x^{i_{1}i_{2}i_{3}i_{4}}y^{i_{1}i_{2}i_{3}i_{4}}$ as $e_{j}$ with $e_{j}\in C_{j}$, $1\leqslant j \leqslant 31$, we can construct a cycle of length $l$ in $BSQ_{n}$ in the following way:
$$C_{l} = (\hat{C}-x^{0000}y^{0000})\cup (C'-x^{0000}y^{0000}-\cup_{j=1}^{31}{e_{j}})\cup (\cup_{j=1}^{31}{(C_{j}-e_{j})}).$$

Hence, $BSQ_{n}$ is vertex-bipancyclic for all $n\geqslant 2$.} \qed

\section{Conclusions}

The simplified shuffle-cube and the balanced shuffle-cube, proposed by  L\"{u} et al. \cite{ref8}, have some properties which superior to the shuffle-cube. In this paper, we prove that the $n$-dimensional simplified shuffle-cube $SSQ_{n}$ is vertex-pancyclic for $n\geqslant 6$ and the $n$-dimensional balanced shuffle-cube $BSQ_{n}$ is vertex-bipancyclic for $n\geqslant 2$. The (bi)pancyclicity of shuffle-cube variants with faulty elements is of interest and should be further studied.
\vskip 0.3 in

\noindent{\bf\large Declarations}

\vskip 0.3 in
\noindent{\bf\normalsize Ethical Approval}
Not applicable.

\vskip 0.3 in

\noindent{\bf\normalsize Funding} 
This research was partially supported by the National
Natural Science Foundation of China (No. 11801061).

\vskip 0.3 in

\noindent{\bf\normalsize Availability of data and materials}
Not applicable.


\vskip 0.05 in

\begin{thebibliography}{99}

\bibitem{TW} S. Abraham, K. Padmanabhan, The twisted cube topology for multiprocessors: a study in network asymmetry, J. Parallel Distrib. Comput. 13 (1991) 104–110.

\bibitem{ref7} S. Anantapantula, C. Melekian, E. Cheng, Matching preclusion for the shuffle-cubes, Parallel Process. Lett. 28 (2018) 1850012.

\bibitem{AM} W.C. Athas, C.L. Seitz, Multicomputers: message-passing concurrent computers, Comput. 21 (1988) 9–24.

\bibitem{ref9} J.A. Bondy, Pancyclic graphs, J. Combin. Theory Ser. B. 11 (1) (1971) 80–84.

\bibitem{AU} S.A. Choudum, V. Sunitha, Augmented cubes, Netw. 40 (2002) 71–84.

\bibitem{ref2} T. Ding, P. Li, M. Xu, The component (edge) connectivity of shuffle-cubes, Theor. Comput. Sci. 835 (2020) 108–119.

\bibitem{TH} T.H. Dunican, Performance of the Intel iPSC/860 and ncube 6400 hypercubes, Parallel Comput. 17 (1991) 1285–1302.

\bibitem{CR} K. Efe, A variation on the hypercube with lower diameter, IEEE Trans. Comput. 40 (1991) 1312–1316.

\bibitem{FO} A. El-Amawy, S. Latifi, Properties and performance of folded hypercubes, IEEE Trans. Parallel Distrib. Syst. 2 (1991) 31–42.

\bibitem{JP} J.P. Hayes, T.N. Mudge, Hypercube supercomputers, in Proc. IEEE. 77 (1989) 1829–1841.

\bibitem{FI} W.-J. Hsu, Fibonacci cubes–a new interconnection topology, IEEE Trans. Parallel Distrib. Syst. 4 (1993) 3–12.

\bibitem{ref1} T.-K. Li, Jimmy J.M. Tan, L.-H. Hsu, T.-Y. Sung, The shuffle-cubes and their generalization, Inform. Process. Lett. 77 (2001) 35–41.

\bibitem{ref8} H. L\"{u}, K. Deng, X. Yang, Symmetric properties and two variants of
shuffle-cubes, arXiv: 2110.13645, 2021.

\bibitem{LU} E. Munarini, C.P. Cippo, N.Z. Salvi, On the Lucas cubes, Fibonacci Quart. 39 (2001) 12–21.

\bibitem{ref3} X.-W. Qin, R.-X. Hao, Reliability analysis based on the dual-CIST in shuffle-cubes, Appl. Math. Comput. 397 (2021) 125900.



\bibitem{ref10} B. Randerath, I. Schiermeyer, M. Tewes, L. Volkmann, Vertex pancyclic graphs, Discrete Appl. Math. 120 (1–3) (2002) 219–237.


\bibitem{QN} Y. Saad, M.H. Schultz, Topological properties of hypercubes, IEEE Trans. Comput. 37 (7) (1988) 867–872.

\bibitem{CO} C.L. Seitz, The Cosmic Cube, Comm. ACM. 28 (1985) 22–33.


\bibitem{BA} J. Wu, K. Huang, The balanced hypercube: a cube-based system for fault-tolerant applications, IEEE Trans. Comput. 46 (1997) 484–490.


\bibitem{AS} J.-M. Xu, M. Ma, A survey on cycle and path embedding in some networks, Front. Math. China 4 (2) (2009) 217–252.

\bibitem{ref4} J.-M. Xu, M. Xu, Q. Zhu, The super connectivity of shuffle-cubes, Inf. Process. Lett. 96 (2005) 123–127.




\end{thebibliography}
\end{document}